# Enhancing System Flexibility through Corrective Demand Response in Security-Constrained Unit Commitment

Arun Venkatesh Ramesh, *Student Member, IEEE* and Xingpeng Li, *Member, IEEE*
Department of Electrical and Computer Engineering
University of Houston
Houston, TX, USA
(Email: aramesh4@uh.edu, xli82@uh.edu)

*Abstract*—Currently, system operators implement demand response by dispatching controllable loads for economic reasons in day-ahead scheduling. Particularly, demand shifting from peak hours when the cost of electricity is higher to non-peak hours to maintain system reliability by flattening the load profile. However, the system flexibility and economic benefits of such action in post-contingency scenarios are not explicitly considered in short-term operations. Hence, this paper highlights the benefits of demand response as a corrective action for potential post-contingency emergencies in day-ahead scheduling. A security-constrained unit commitment (SCUC) model which considers the flexibility offered through corrective demand response (CDR) to maintain system reliability when a line or generator outage occurs is proposed. The proposed model was tested on IEEE 24-bus system where simulation results point to significant total cost savings in daily operations. Moreover, the results point to better long-term reliability of generators along with the ability to use existing system flexibility and serve higher critical demands in base-case when CDR is implemented.

## Nomenclature

*Sets/Indices:*

| | |
|---|---|
| $g$ | Generator index. |
| $G$ | Set of generators. |
| $g(n)$ | Set of generators connecting bus $n$. |
| $k$ | Line index. |
| $K$ | Set of all transmission element. |
| $\delta^+(n)$ | Set of lines with bus $n$ as receiving bus. |
| $\delta^-(n)$ | Set of lines with bus $n$ as sending bus. |
| $t$ | Time period index. |
| $T$ | Set of Time intervals. |
| $n$ | Bus index. |
| $N$ | Set of all buses. |
| $N(g)$ | Bus location of generator $g$. |
| $c$ | Contingent element index. |
| $C$ | Set of contingencies. |
| $G_c$ | Set of all generator contingencies. |
| $K_c$ | Set of all non-radial line contingencies. |

*Parameters:*

| | |
|---|---|
| $UT_g$ | Minimum up time for generator $g$. |
| $DT_g$ | Minimum down time for generator $g$. |
| $c_g$ | Linear cost for generator $g$. |
| $c_g^{NL}$ | No-load cost for generator $g$. |
| $c_g^{SU}$ | Start-up cost for generator $g$. |
| $\pi_c$ | Probability of contingency $c$. |
| $c_n^{Ctg}$ | Cost of CDR at bus $n$. |
| $P_g^{min}$ | Minimum capacity of generator $g$. |
| $P_g^{max}$ | Maximum capacity of generator $g$. |
| $R_g^{hr}$ | Regular hourly ramping limit of generator $g$. |
| $R_g^{SU}$ | Start-up ramping limit of generator $g$. |
| $R_g^{SD}$ | Shut-down ramping limit of generator $g$. |
| $R_g^{10}$ | 10-minute outage ramping limit of generator $g$. |
| $P_k^{max}$ | Long-term thermal line limit for line $k$. |
| $b_k$ | Susceptance of line $k$. |
| $P_k^{emax}$ | Emergency thermal line limit for line $k$. |

*Variables:*

| | |
|---|---|
| $P_{g,t}$ | Output of generator $g$ in time period $t$. |
| $u_{g,t}$ | Commitment status of generator $g$ in time period $t$. |
| $v_{g,t}$ | Start-up variable of generator $g$ in time period $t$. |
| $r_{g,t}$ | Reserve from generator $g$ in time period $t$. |
| $P_{k,t}$ | Lineflow of line $k$ in time period $t$. |
| $\theta_{ref,t}$ | Phase angle of reference bus in time period $t$. |
| $\theta_{n,t}$ | Phase angle of bus $n$ in time period $t$. |
| $\theta_{m,t}$ | Phase angle of bus $m$ in time period $t$. |
| $d_{n,t}$ | Demand of bus $n$ in time period $t$. |
| $CDR_{n,c,t}$ | Corrective demand response action at bus $n$ in period $t$ for contingency $c$ |
| $P_{g,c,t}$ | Output of generator $g$ in period $t$ for contingency $c$ |
| $P_{k,c,t}$ | Flow in line $k$ in period $t$ after outage of equipment $c$. |
| $\theta_{m,c,t}$ | Phase angle of bus $m$ in period $t$ for contingency $c$. |
| $\theta_{n,c,t}$ | Phase angle of bus $n$ in period $t$ for contingency $c$. |

## I. Introduction

With the advent of smart grid technologies, we have brought about intelligent energy management systems to optimally and reliably operate the grid. Traditionally, grid was operated in a top-down framework where the flexibility requirement of the power system is met with committing additional generators to meet the demand and the reliability requirement of the network. But the technology to sense and control signals with two-way communications has brought increased participation from demand side in energy markets [1]. The system operators can also determine and send signals to not only redispatch generators but also adjust controllable loads.

Utilities now offer several price-based or incentive-based programs for altering demand patterns through demand side management (DSM) that reduces the cost of the electricity the customer pays [2]. DSM not only lowers cost but also enhances reliability and provides self-healing capabilities for the power grid through demand response (DR) [3]. In particular, demand response through direct load control (DLC) enables grid operators to send signals to reduce non-critical loads directly. However, most DLC actions are implemented as part of the distribution network by utilities to shift non-critical loads from peak hours experiencing high demands to non-peak hours [2].

In day-ahead scheduling, system operators use security-constrained unit commitment (SCUC) to obtain the optimal commitment status and dispatch signals for generators to meet forecasted bulk hourly loads [4]. As per Federal Energy Regulatory Commission (FERC), the SCUC solution requires to be *N-1* reliable where the system is capable of handling frequent line or generator outages individually [5].

Here, system operators often utilize preventive and corrective control actions to maintain system reliability [4], [6]-[7]. Mostly, DR through controllable loads are considered as a preventive action. DR benefits the system by moving non-critical deferrable loads from peak hours to non-peak hours which increases the system flexibility and demand side market participation [1]-[2]. Though there are emergency DR plans that several system operators implement, they are solely based on supply and demand balance for frequency regulation and typically this is factored in the SCED process where operators dispatch the participating DR resources [8].

Mostly, independent system operators (ISOs) compensate DLC by locational marginal prices (LMP) when utilized in real-time operations and fixed-price based schemes contracted in long-term capacity markets to shift demand [9]-[13]. By studying the ISO reports for outage statistics for the network with about 350 dispatchable generators and 9,000 miles of high voltage transmission lines, it reported that a low number of cases, 833, for the unplanned outages related to line in 2019 and unplanned outages generator are very rare [14],[15]. Here, the untapped potential of DLC considering system flexibility or capacity release is not completely considered in short-term operations. Mostly, DR through dispatchable DLC are used for economic reasons and reliability reasons in base-case; it is only used for reliability purpose in post-contingency scenario.

Most utilization of DLC is preventive in nature as seen in [16]-[17]. In [18], the DR actions were automated to respond to real-time dispatch schedule to provide additional reserves as an ancillary service to the grid. The DLC algorithms consider economic benefits in [19] where optimal control schedules are determined by nodal aggregators by optimizing the load profile. Many preventive DR actions also consider the reliability of power system but with limited contingencies being considered as seen in [20]. However, for economic benefits of DR actions, they need to be implemented in operational optimization problem such as SCUC or security-constrained economic dispatch (SCED) and only [4], [19] and [21] include the power network constraints. Here, all the research address DR as a preventive action only, which may substantially affect customer comfort level and are susceptible to cyber security in real-time operations [22].

Currently, the use of corrective demand response (CDR) in response to contingency is never considered in the SCUC process. Similar to network reconfigurations as a corrective action [7], [23]-[25], CDR can also increase the solution quality by reducing costs when co-optimized with SCUC. Not only that, CDR can provide additional system flexibility by shedding some non-critical load under contingency and allowing the committed units to ramp-up or ramp-down to meet the system requirements in post-contingency scenarios rather than committing additional units. Note that under CDR schemes, non-critical load shedding will occur only when the associated contingency actually occurs, which is a low probably event.

Though there exists load shedding, especially during low supply scenarios due to congestion-induced or fault-induced events, this is however implemented for reliability and not for economic benefits in [6]. In [26], DR is considered as a corrective action in post-fault condition to release network capacity but not for economic benefits. Similarly, in [27], a few DLC operations are considered in a post-contingency scenario to obtain additional system flexibility to enhance system reliability under an emergency but with limited economic considerations. In [26]-[27], the operational and network constraints are not completely considered or are only focused on distribution networks. Therefore, the studies on benefits of utilizing DLC in day-ahead operations at the transmission level for post-contingency actions are limited.

In this paper, Section I addresses the research gap and Section II then formulates the proposed model. Section III describes the test system that is used to validate the proposed model whereas Section IV presents the experimental results achieved along with detailed analysis. Finally, Section V concludes the paper and Section VI provides the future scope of the proposed research.

## II. PROPOSED MODEL

The proposed mixed-integer linear programming model of SCUC-CDR minimizes the total cost (1) consisting of operational cost of generators and the penalty cost for the CDR action while satisfying the base-case and reliability constraints. The penalty cost does not affect the total cost of SCUC formulation since it does not perform any CDR actions. The base-case constraints in (2)-(15) combines both generation and power flow constraints. The generator physical constraints consist of, (i) maximum and minimum limits on generations (3) and (2); (ii) reserve requirements are considered in (4) and (5); (iii) the ramping limits for hourly time-period is represented in (6) and (7); (iv) minimum-down time before a generator can be started-up and the minimum-up time before a generator can be shutdown are depicted in (9) and (8), respectively; (v) the generator start-up is defined in (10) and the indicating variables for generator start-up and commitment status are binary and are represented by (11). The physical power flow constraints consist of (i) nodal power balance which meets supply and demand, (12); (ii) line flow limits enforced in (13); (iii) line flow equations, (14); (iv) the system reference node is defined in (15).

*Objective:*

$$\text{Min} \sum_g \sum_t (c_g P_{g,t} + c_g^{NL} u_{g,t} + c_g^{SU} v_{g,t}) + \sum_{n,c,t} (\pi_c * c_n^{Ctg} * CDR_{n,c,t}) \quad (1)$$

s.t.:

*Base case modeling of generation:*

$$P_g^{min} u_{g,t} \leq P_{g,t}, \forall g, t \quad (2)$$
$$P_{g,t} + r_{g,t} \leq P_g^{max} u_{g,t}, \forall g, t \quad (3)$$
$$0 \leq r_{g,t} \leq R_g^{10} u_{g,t}, \forall g, t \quad (4)$$
$$\sum_{q \in G} r_{q,t} \geq P_{g,t} + r_{g,t}, \forall g, t \quad (5)$$
$$P_{g,t} - P_{g,t-1} \leq R_g^{hr} u_{g,t-1} + R_g^{SU} v_{g,t}, \forall g, t \quad (6)$$
$$P_{g,t-1} - P_{g,t} \leq R_g^{hr} u_{g,t} + R_g^{SD}(v_{g,t} - u_{g,t} + u_{g,t-1}), \forall g, t \quad (7)$$
$$\sum_{q=t-UT_g+1}^{t} v_{g,q} \leq u_{g,t}, \forall g, t \geq UT_g \quad (8)$$
$$\sum_{q=t+1}^{t+DT_g} v_{g,q} \leq 1 - u_{g,t}, \forall g, t \leq T - DT_g \quad (9)$$
$$v_{g,t} \geq u_{g,t} - u_{g,t-1}, \forall g, t \quad (10)$$
$$v_{g,t}, u_{g,t} \in \{0,1\}, \forall g, t \quad (11)$$

*Base case modeling of power flow:*

$$\sum_{g \in g(n)} P_{g,t} + \sum_{k \in \delta^+(n)} P_{k,t} - \sum_{k \in \delta^-(n)} P_{k,t} = d_{n,t}, \forall n, t \quad (12)$$
$$-P_k^{max} \leq P_{k,t} \leq P_k^{max}, \forall k, t \quad (13)$$
$$P_{k,t} - b_k(\theta_{n,t} - \theta_{m,t}) = 0, \forall k, t \quad (14)$$
$$\theta_{ref,t} = 0 \;\forall t \quad (15)$$

The reliability constraints when a contingency occurs are modelled in (16)-(22), typically with a contingency response time frame of 10 minutes. Those consist of both the generation constraints, (16)-(19), and power flow constraints, (20)-(22) in post-contingency scenarios. As per the post-contingency generator constraints, 10-minute ramp up/down limits of generators are factored in (16)-(17) while adhering to minimum and maximum physical limits of generator in (18)-(19). The post-contingency line flows are determined in (20) while adhering to emergency line limits in (21). The nodal power balance in post-contingency scenario is enforced in (22) whereas the nodal power balance in post-contingency scenario considering the corrective demand shedding is enforced in (23). Finally, (24) defines that the corrective demand shedding at each participating node is assumed to be capped at 30% of the nodal demand.

*Post-contingency 10-minute ramping restriction on generation and modeling of contingencies:*

$$P_{g,t} - P_{g,c,t} \leq R_g^{10} u_{g,t}, \forall g, c \in C, t \quad (16)$$
$$P_{g,c,t} - P_{g,t} \leq R_g^{10} u_{g,t}, \forall g, c \in C, t \quad (17)$$
$$P_g^{min} u_{g,t} \leq P_{g,c,t}, \forall g, c \in C, t \quad (18)$$
$$P_{g,c,t} \leq P_g^{max} u_{g,t}, \forall g, c \in C, t \quad (19)$$

*Post-contingency modeling of power flow:*

$$P_{k,c,t} - b_k(\theta_{n,c,t} - \theta_{m,c,t}) = 0, \forall k, c \in C, t \quad (20)$$
$$-P_k^{emax} \leq P_{k,c,t} \leq P_k^{emax}, \forall k, c \in C, t \quad (21)$$
$$\sum_{g \in g(n)} P_{g,c,t} + \sum_{k \in \delta^+(n)} P_{k,c,t} - \sum_{k \in \delta^-(n)} P_{k,c,t} = d_{n,t} \quad (22)$$
$$\sum_{g \in g(n)} P_{g,c,t} + \sum_{k \in \delta^+(n)} P_{k,c,t} - \sum_{k \in \delta^-(n)} P_{k,c,t} = d_{n,t} - CDR_{n,c,t} \quad (23)$$
$$CDR_{n,c,t} \leq 0.3 * d_{n,t} \quad (24)$$

Since generator contingencies are very rare compared to line contingencies, different models are formulated for both SCUC and SCUC-CDR. Mainly, T-SCUC and TG-SCUC models the SCUC with only transmission contingencies and with both transmission and generator constraints, respectively. Whereas, the proposed models for SCUC with corrective actions using DLC for post-contingency constraints are namely, T-SCUC-CDR and TG-SCUC-CDR, with transmission contingencies only and both transmission and generator contingencies, respectively. Here, T-SCUC and TG-SCUC are defined by (1)-(22). T-SCUC-CDR and TG-SCUC-CDR are modelled through (1)-(21) and (23)-(24). The difference in T and TG models are captured through input set of contingencies, $C$, where $C \in K_c$ for transmission contingencies only and $C \in G_c \cup K_c$ when both transmission and generator contingencies are modelled. Based on the above constraints, the proposed models are consolidated in Table I.

TABLE I. PROPOSED MODEL CONSTRAINTS

| Model | T-SCUC | TG-SCUC | T-SCUC-CDR | TG-SCUC-CDR |
|---|---|---|---|---|
| Objective | (1) | (1) | (1) | (1) |
| Constraints | (2)-(22) | (2)-(22) | (2)-(21), (23)-(24) | (2)-(21), (23)-(24) |
| $C$ | $K_c$ | $G_c \cup K_c$ | $K_c$ | $G_c \cup K_c$ |

### III. TEST CASE DESCRIPTION

The proposed SCUC-CDR was validated against traditional SCUC on the IEEE 24-bus system [28]. This test system contains 24 buses, 33 generators and 38 lines as shown in Fig. 1. The total generation capacity from generators is 3,393 MW and the system peak load is 2,281 MW, [23].

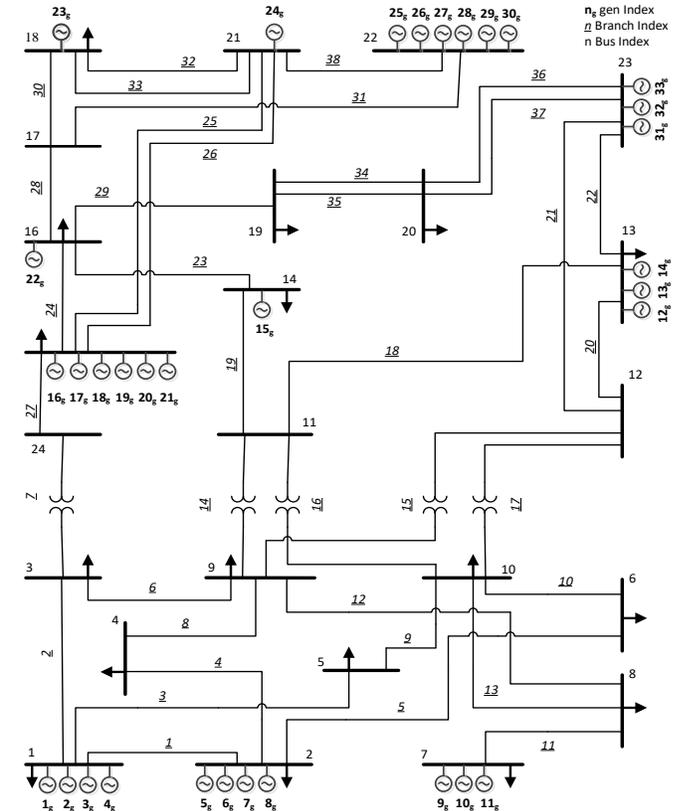

Fig. 1 IEEE 24-bus system network diagram [23].

In this paper, the probability of any transmission outage in the network is the inverse of total number of non-radial lines. Similarly, the probability of any generator outage in the system is the inverse of total number of generators. Note that the probability of the above outages are quite rare and the realistic contingency probability can be modelled through long-term outage statistics.

## IV. RESULTS

The mathematical model was implemented using AMPL and solved using Gurobi solver for a 24-hour (Day-Ahead) load period for the test system described in Section III.

### A. Total economic benefits of CDR:

The difference in overall cost of T-SCUC-CDR and TG-SCUC-CDR against the T-SCUC and TG-SCUC demonstrates a cost saving of $9,825 and $14,996, respectively when CDR is introduced. This is due to the flexibility offered by CDR actions which provides a more economical commitment status and fewer generator start-ups to handle the same demand. Also, the more constrained problem which considers both line and generator contingencies, TG-SCUC and TG-SCUC-CDR, results in a higher cost saving than the respective models, T-SCUC and T-SCUC-CDR, that only consider transmission outages. However, it can be noted that generator outages are very infrequent compared to line outages.

The proposed model considering only transmission contingencies, T-SCUC-CDR, results in a total curtailment of 25.5 MW over 24 hours as CDR action which is same in the case of TG-SCUC-CDR. It can be noted that the cumulative CDR action for line outages is only 0.01% of the peak system load and it brings about significant total operational cost reduction. It was also observed that key system line, line 7 or line 27, outage required 8.4 MW CDR action at bus 14; and line 8 outage resulted with 2.18 MW CDR action at bus 6 at various time periods. Therefore, only few critical outages required CDR to satisfy system requirements. TG-SCUC-CDR, the total amount of CDR action is much higher, 370 MW over 24 hours for generator outages and 25.5 MW over 24 hours for line outages. Similar to CDR due to line outages, only a few large generator outages, generator 23 and generator 24, utilize CDR to maintain system reliability.

Both TG-SCUC-CDR and T-SCUC-CDR models benefit by significantly faster solve time when compared to TG-SCUC and T-SCUC, respectively. In particular, TG-SCUC-CDR is 20% faster than TG-SCUC and T-SCUC-CDR is 48% faster than T-SCUC. This is because, the introduction of CDR results in a relaxed problem with increased feasible set of solutions.

TABLE II. OPERATIONAL COST AND POST-CONTINGENCY DEMAND CURTAILED

|  | TG-SCUC | TG-SCUC-CDR | T-SCUC | T-SCUC-CDR |
|---|---|---|---|---|
| Cost ($) | 685,670 | 670,674 | 677,851 | 668,026 |
| MIPGAP | 0.0095 | 0.0045 | 0.0085 | 0.0015 |
| Solve time (s) | 237 | 191 | 111 | 57 |
| $\sum_{n,c \in K_{c,t}} CDR_{n,c,t}$ | NA | 25.46 | NA | 25.46 |
| $\sum_{n,c \in G_{c,t}} CDR_{n,c,t}$ | NA | 369.81 | NA | NA |

### B. CDR penalty cost sensitivity:

Since the CDR actions result in curtailment of non-critical loads, a penalty cost for such actions are introduced in the objective cost in (1). Moreover, the occurrences of either transmission or generator outages are low which is modelled by $\pi_c$. The product of the probability, $\pi_c$, and cost of CDR, $c_n^{Ctg}$, represents the penalty cost in the system. The system was studied with varying cost of CDR from 0 $/MWh to 40,000 $/MWh. This is represented in two graphical forms, a low penalty cost sensitivity to CDR actions, Fig. 2 and a high penalty cost sensitivity to CDR actions, Fig. 3. The system shows inverse relations, that is the cumulative amount CDR actions decreases as cost of CDR increases.

In the low penalty cost sensitivity, at 0 $/MWh, there is no control to limit CDR actions and hence it results in cumulative curtailment of 328,925 MW for line outages and 341,379 MW for generator outages (not represented in the scale of graph). At 1 $/MWh, we notice significant reduction to the CDR action with 25.46 MW for line outages and 369.8 MW for generator outages. It was noted that at low cost of CDR, the CDR actions due to generator outages are more sensitive whereas due to line outages are constant. Here, the total cost of the system changed marginally to increasing penalty cost with the anomaly at 1 $/MWh can be explained by the associated higher relative gap in solution.

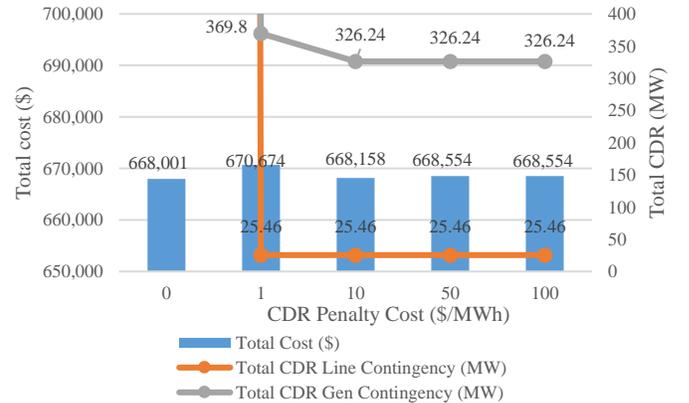

Fig. 2 Low penalty cost sensitivity study for TG-SCUC-CDR.

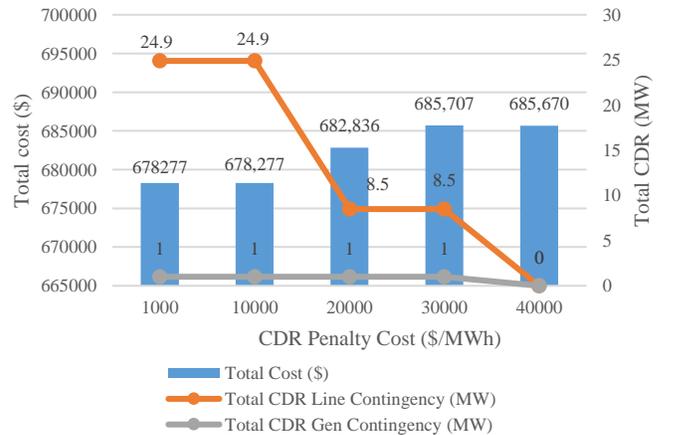

Fig. 3 High penalty cost sensitivity study for TG-SCUC-CDR.

A high penalty cost sensitivity study was conducted to identify when the total CDR actions result 0 MW. Here, it was noted that only at very high cost of CDR, $40,000, the system

does not implement CDR for both line outages and generator outages. Therefore, the total cost of the system for TG-SCUC-CDR is same as the TG-SCUC. Also, it was noted that at high cost of CDR, the CDR actions due to generator outages and line outages are very sensitive. The cumulative CDR actions due to line outages dropped steeply from 25.46 MW to 8.5 MW at a penalty cost of $540 whereas due to generator outages dropped steeply at a penalty cost of $30 from 326.24 MW to 1 MW.

*C. System flexibility:*

Five scenarios were considered: two low-load scenarios (80%, 90%), a base-load scenario (100%) and two high-load scenarios (110%, 120%). The load profile was varied using a percentage multiplied to the nodal load. Fig. 4 shows the total cost for various methods under different load profiles and respective cumulative CDR actions for generator and line outages.

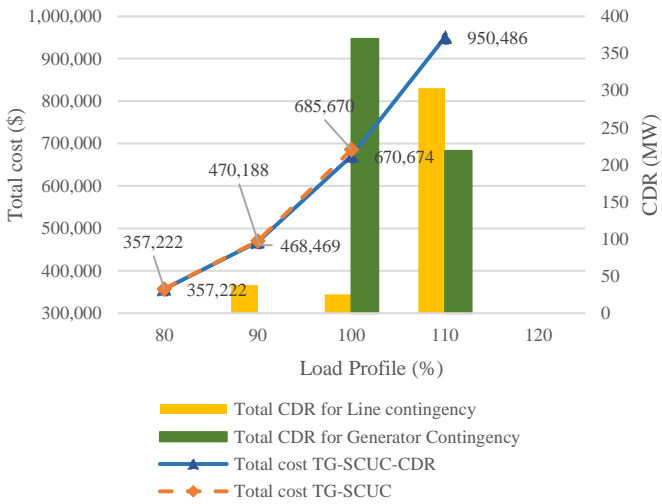

Fig. 4 Total system cost and cumulative DR shifted for different scenarios.

CDR is never implemented for the very low-load scenario (80%) since the base-case network loading level is low and post-contingency networks are not congested. This implies both TG-SCUC and TG-SCUC-CDR obtain the same total cost. At very high load scenario, 120%, TG-SCUC and TG-SCUC-CDR are infeasible.

As the system is loaded, CDR actions are observed in load scenarios of 90%-110% along with economic benefits of total cost reduction by cheaper generator dispatch schedule. Here, at 90% load scenario only line outages resulted in CDR actions whereas CDR actions in base-load (100%) and high-load scenarios (110%) are resulted from both line and generator outages. At base-load scenario (100%), the cumulative CDR actions due to line outages are lower compared to cumulative CDR actions due to generator outages.

However, at high-load scenario (110%), TG-SCUC is infeasible whereas TG-SCUC-CDR provides a feasible solution by utilizing the flexibility in the system associated with CDR. This implies that utilizing CDR is beneficial in serving higher critical loads compared to traditional SCUC which does not implement any corrective actions.

*D. Market analysis:*

Table III shows the market results for base-load profile (100%) which compare the load payment, generator revenue and average nodal LMP for various scheduling models when CDR is utilized (T-SCUC-CDR, TG-SCUC-CDR) and when CDR is not utilized (T-SCUC, TG-SCUC). Overall, it is observed that with CDR the average nodal LMP, load payment and generator revenue are higher, which is counter-intuitive since TG-SCUC-CDR or T-SCUC-CDR results in a lower total operation cost solution. Also, the difference in average nodal LMP is more evident when both line and generator outages are considered compared to only line outages.

The higher LMP in TG-SCUC-CDR compared to TG-SCUC can be explained using the generator commitment solution in Table IV. Since the market results are calculated with LMP, it is expected to have higher load payment and generator revenue due to higher average nodal LMP. However, the commitment solution, 14 for TG-SCUC and 6 for TG-SCUC-CDR after period 1, favors long-term reliability of generators through infrequent generator start-ups with flexibility obtained through CDR. Here, all units are OFF before period 1; and in period 1, both TG-SCUC and TG-SCUC-CDR commit 18 units. However, there are more uncommitted units which are always OFF in TG-SCUC-CDR compared to TG-SCUC. Traditionally, the flexibility in the system is obtained by committing extra units as seen in TG-SCUC where a total of 474 committed generator-hours over 24 hours was noticed whereas in TG-SCUC-CDR, it was bettered efficiently to 460 committed generator-hours over 24 hours.

The nodal LMP is higher in the case of TG-SCUC-CDR due to: (i) fewer generators are committed (ii) marginal units are more expensive (iii) the cheaper generators (always ON) capacity are completely utilized in the base-case dispatch solution. There are also more expensive units that are in 'always OFF' condition in TG-SCUC-CDR, which points to reduced generator start-ups.

TABLE III. MARKET RESULTS FOR IEEE-24 BUS SYSTEM

|  | TG-SCUC | TG-SCUC-CDR | T-SCUC | T-SCUC-CDR |
|---|---|---|---|---|
| Load payment($) | 1,289,650 | 1,698,060 | 1,070,370 | 1,071,270 |
| Gen revenue ($) | 707,594 | 1,073,830 | 1,683,830 | 1,695,420 |
| Avg LMP ($) | 23.7 | 31.74 | 31.51 | 31.67 |

TABLE IV. GENERATOR COMMITMENT STATUS

|  | TG-SCUC | TG-SCUC-CDR |
|---|---|---|
| Always ON | 3,7,8,21-33 | 4,7-8,21-33 |
| Always OFF | 12-15 | 1-2, 5-6,15-16,19-20 |
| Marginal Units | 1-6,9-11,16-20 | 3,9-11,17-18 |
| Total start-ups $t > 1$ | 14 | 6 |
| Total start-ups $t = 1$ | 18 | 18 |
| Total commitment | 474 | 460 |

## V. CONCLUSION

The use of demand response as a corrective action to system contingencies was proposed as TG-SCUC-CDR model and studied. Mainly, TG-SCUC-CDR and T-SCUC-CDR results in lower operational costs by reducing generator start-ups and fewer generators committed for the same load profile compared to TG-SCUC and T-SCUC, respectively. In particular, the sensitivity of such CDR actions which provide significant economic benefits were studied with respect to penalty cost and load profile variation. The results indicate that given a high demand profile, SCUC is infeasible whereas SCUC-CDR is feasible as the system uses the available system flexibility.

Also, the sensitivity to penalty cost shows that CDR actions with even small amounts such as 2 MW can result in substantial economic benefits. It was also noted that the CDR actions due to generator outages are more sensitive to variation in penalty costs.

The market analysis resulted in counter-intuitive results as the average nodal LMP were higher for TG-SCUC-CDR and T-SCUC-CDR compared to TG-SCUC and T-SCUC, respectively. However, this was explained by the additional units committed and the capacity of cheaper generators were not completely exhausted at the cost of expensive units running at no-load or low capacities in the case of TG-SCUC and T-SCUC.

## VI. Future work

The future work can consider the interaction of CDR action which provides additional system flexibility with the variability associated with renewable energy such as wind and solar, which can be of importance. Especially, the proposed model with corrective DR actions can be implemented in stochastic model to consider the various weather scenarios like the use of corrective network reconfiguration (CNR) action in a stochastic model [29].

The interaction of DR/CDR actions with other technologies such as energy storage and network topology modification efficiently is another topic to be researched.